\documentclass[10pt,a4paper,draft]{article}

\usepackage[english]{babel}
\usepackage[intlimits]{amsmath}
\usepackage[]{amsfonts}
\usepackage[]{amssymb}
\usepackage[]{a4wide}
\usepackage[]{amsthm}

\newtheorem{thm}{Theorem}
\newtheorem{lem}[thm]{Lemma}
\theoremstyle{definition}\newtheorem{rem}[thm]{Remark}
\theoremstyle{definition}\newtheorem{exa}[thm]{Example}

\newcommand{\dom}{\mathrm{dom}\:}
\newcommand{\ran}{\mathrm{ran}\:}
\newcommand{\fdom}[1]{\mathrm{D_\ast}\:(#1)}
\newcommand{\adj}{^\ast}
\newcommand{\aadj}{^{\ast\ast}}
\newcommand{\nam}{\textsc}
\newcommand{\im}{\mathsf{M}}
\newcommand{\hi}{\mathfrak{H}}
\newcommand{\kr}{{\mathsf{K}}}
\newcommand{\fr}{{\mathsf{F}}}
\newcommand{\fv}{\rightarrow}
\newcommand{\ls}{(}
\newcommand{\rs}{)}
\newcommand{\la}{\langle}
\newcommand{\ra}{\rangle}
\newcommand{\fsum}{\stackrel{.}{+}}
\newcommand{\sq}{^{\frac{1}{2}}}
\newcommand{\nm}{\|}
\newcommand{\psup}{\sup}
\newcommand{\negs}{\!\!\!\!\!}
\newcommand{\bo}{\mathcal{B}}
\newcommand{\rest}{\,\rule[-1ex]{0.1ex}{2.5ex}\,}
\newcommand{\spa}{\enskip}
\newcommand{\nul}{\mathbf{0}}
\newcommand{\cc}{\mathbb{C}}

\newenvironment{ix}{\begin{array}{c}}{\end{array}}

\begin{document}


\title{Commutation properties of the form sum of positive, symmetric
operators}

\author{
B\'alint Farkas\\
\small\noindent\nam{Department of Applied Analysis, E\"otv\"os Lor\'and University}\\
\small\noindent {Kecskem\'eti u. 10-12, 1053 Budapest, Hungary}\\
\small\noindent e-mail: \texttt{fbalint@cs.elte.hu}
\\
\tiny{and}\\
M\'at\'e Matolcsi\\
\small\noindent\nam{Department of Applied Analysis, E\"otv\"os Lor\'and University}\\
\small\noindent {Kecskem\'eti u. 10-12, 1053 Budapest, Hungary}\\
\small\noindent e-mail: \texttt{matomate@cs.elte.hu}
}
\maketitle

\begin{abstract}A new construction for the form sum of positive, selfadjoint operators is given in this paper. The situation is a bit more general, because our aim is to add positive, symmetric operators. With the help of the used method, some commutation properties of the form sum extension are observed.
\end{abstract}
\section{Introduction}
Given two positive, selfadjoint operators $A$ and $B$ in the Hilbert space $\hi$, we may form the operator sum $A+B$ on $\dom{A}\cap\dom{B}$. However, the intersection of the domains may be zero-dimensional, and in general nothing can assure us that the sum will be a selfadjoint operator. The so called form sum extension handles this problem if $\dom{A\sq}\cap\dom{B\sq}$ is dense in $\hi$. Define $q_A(x)=\ls A\sq x,A\sq x\rs$ and $q_B(x)=\ls B\sq x,B\sq x\rs$ two closed forms; their sum $q_A+q_B$ is a closed form on $\dom{A\sq}\cap\dom{B\sq}$, therefore the representation theorem provides us a selfadjoint operator $C$, such that $C$ and $A+B$ coincide on $\dom{A}\cap\dom{B}$ \cite{faris}. The usual notation for the form sum of $A$ and $B$ is $A\fsum B$. In Section 2, we give a new construction of the form sum of positive, symmetric operators. Section 3 deals with commutation properties of this extension, i.e. how our extension method can preserve commutation with bounded operators. In the last section we give some examples concerning the form sum extension, and describe the relation between other extensions of operator sums.

We use the following notations, and refer the reader to \cite{sebprok2}, \cite{sebprok1} and \cite{sebstoch}. Throughout $a,b$ will denote positive, symmetric operators in the Hilbert space $\hi$, with not neccesarily dense domains.
$\fdom a$ will denote the so-called form domain of $a$, i.e. 
$$\fdom{a}=\{y\in\hi:\exists m_y\spa |\ls ax,y\rs|^2\leq m_y\ls ax,x\rs,\forall x\in\dom{a}\}.$$
We remark, if $a$ is positive, selfadjoint, then $\fdom{a}=\dom{a\sq}$.
The \nam{Krein-von Neumann} and \nam{Friedrichs} extensions of $a$ will be denoted by $a_\kr$ and $a_\fr$ respectively (provided they exist).
We recall the basic notions now. If $\fdom{a}$ is dense in $\hi$, then  $\la ax,ay\ra=\ls ax,y\rs$ is an inner product on $\ran{a}$. Let $\hi_a$ denote the completion of the pre-Hilbert space $\ran a$ with the above inner product. Define $J_a:\hi_a\rightarrow\hi$ and $Q_a:\hi\rightarrow\hi_a$ by 
\begin{gather*}
\dom{J_a}=\ran{a},\quad J_aax=ax \mbox{ for all }ax\in\ran{a}\\
\dom{Q_a}=\dom{a},\quad Q_ax=ax \mbox{ for all }x\in\dom{a}\\
\end{gather*}
Now, if $\fdom{a}$ is dense, then $J_a\aadj J_a\adj$ is the smallest positive selfadjoint extension of $a$, i.e. the \nam{Krein-von Neumann} extension (see  \cite{sebprok1}, \cite{sebstoch}). Also , the following characterizing properties of $a_\kr$ will be used frequently in this paper
\begin{equation*}\label{eq:-0}
\begin{gathered}
\dom{a_\kr\sq}=\dom{J_a\adj}=\fdom{a}\\
\|a_\kr\sq y\|^2=\|J_a\adj y\|^2=\psup_{\begin{ix}\scriptstyle x\in\dom{a}\\\scriptstyle\la ax,x\ra\leq 1\end{ix}}|(ax,y)|^2\quad\mbox{ for all $y\in\dom{a_\kr\sq}$}
\end{gathered}
\end{equation*}

Provided that $\dom{a}$ is dense, $Q_a\adj Q_a\aadj$ furnishes the largest positive selfadjoint extension, that is the \nam{Friedrichs} extension of $a$ (see \cite{sebprok2}). Note that $\dom{a}\subseteq\fdom{a}$, therefore the denseness of $\dom{a}$ implies the same for $\fdom{a}$.

\section{The form sum extension}
In the following we give a new construction for the sum of two positive, symmetric operators. We show that in case of selfadjont operators this construction supplies the form sum of the operators.  

Let $a$ and $b$ be two positive, symmetric operators, and suppose that $\fdom{a}\cap\fdom{b}$ is dense in $\hi$. Consider the space $\hi_a\oplus \hi_b$, and the operator 
\begin{equation}
\label{eq:0}
J:\hi_a\oplus \hi_b\fv \hi, \mbox{ with } \dom{J}=\ran{a}\oplus\ran{b},\enskip J(ax\oplus by)=ax+by.
\end{equation}
 It is easy to prove that $J\adj$ is densely defined; indeed  $\fdom{a}\cap\fdom{b}=\dom{J\adj}$. To see this, let $x\in\dom{a},y\in\dom{b}$ and $u\in\fdom{a}\cap\fdom{b}$, then
\begin{equation*}
\begin{gathered}
|\ls J(ax\oplus by),u\rs|^2=|\ls ax,u\rs+\ls by,u\rs|^2\leq2|\ls ax,u\rs|^2+2|\ls by,u\rs|^2\leq\\
 2m_u\ls ax,x\rs+2n_u\ls by,y\rs\leq m\la ax\oplus by,ax\oplus by\ra,
\end{gathered}
\end{equation*}
with $m=2\max(m_u,n_u)$. This shows that $u\in\dom{J\adj}$, hence $\fdom{a}\cap\fdom{b}\subseteq\dom{J\adj}$. For the reverse, let $u\in\dom{J\adj}$ and $x\in\dom{a}$, then
$$|\ls ax,u\rs|^2=|\ls J(ax\oplus\nul),u\rs|^2\leq m\la ax\oplus\nul,ax\oplus\nul\ra=m\la ax,ax\ra=m\ls ax,x\rs,$$
with a suitable $m\geq 0$, therefore $u\in\fdom{a}$. Similarly, we obtain that $u\in\fdom{b}$. Thus we have shown that $\fdom{a}\cap\fdom{b}\supseteq\dom{J\adj}$.

We see that $J\aadj$ exists. Now, we calculate $J\adj$ on $\dom{a}\cap\dom{b}$. Let $u\in\dom{a}\cap\dom{b}$ and $x\in\dom{a},y\in\dom{b}$, then
$$\ls J(ax\oplus by),u\rs=\ls ax,u\rs+\ls by,u\rs=\la ax,au\ra+\la by,bu\ra=\la ax\oplus by,au\oplus bu\ra,$$
consequently $J\adj u=au\oplus bu$.

According to the \nam{von Neumann} theorem $J\aadj J\adj$ is positive and selfadjoint. We claim that $J\aadj J\adj$ is an extension of $a+b$. Indeed, let $u\in\dom{a}\cap\dom{b}$, then
$$J\aadj J\adj u=J\aadj (au\oplus bu)=J(au\oplus bu)=au+bu=(a+b)u.$$

In order to prove that our construction is a generalization of the form sum of selfadjoint operators, we need the following lemma on the \nam{Krein-von Neumann} extension (see \cite{sebstoch}, \cite{sebprok2} and \cite{sebprok1}).

\begin{lem}
\label{thm:1}
If $a,b$ are positive, symmetric operators, and $\fdom{a}$ and $\fdom{b}$ are dense in $\hi$, then $\fdom{a\oplus b}$ is dense in $\hi\oplus\hi$ and $$a_\kr\oplus b_\kr=(a \oplus b)_\kr.$$

\begin{proof} First we show that $(a\oplus b)_\kr$ exists. It is enough to prove that $\fdom{a\oplus b}=\dom{(a_\kr\oplus b_\kr)\sq}$ since the latter is dense in $\hi\oplus\hi$.

We observe first that $(a_\kr\oplus b_\kr)\sq=(a_\kr\sq\oplus b_\kr\sq)$, indeed both are positive and selfadjoint with the same square $a_\kr\oplus b_\kr$.

Now, using the definition, we can write:
\begin{equation}
\label{eq:1}
\begin{gathered}
\fdom{a\oplus b}=\\
\{x\oplus y:\exists m_{x,y}\spa|\ls (a\oplus b)(u\oplus v), x\oplus y\rs|^2\leq m_{x,y}\ls  (a\oplus b)(u\oplus v),u\oplus v\rs, \forall u\oplus v\in \dom{a\oplus b}\}\\
=\{x\oplus y:\exists m_{x,y}\spa|\ls au,x\rs+\ls bv,y\rs|^2\leq m_{x,y}(\ls au,u\rs+\ls bv,v\rs) \forall u\oplus v\in \dom{a}\oplus\dom{b}\}.
\end{gathered}
\end{equation}
\begin{equation}
\label{eq:a}
\begin{gathered}
\dom{(a_\kr\oplus b_\kr)\sq}=\dom{(a_\kr\sq\oplus b_\kr\sq)}=\dom{a_\kr\sq}\oplus \dom{b_\kr\sq}=\fdom{a}\oplus\fdom{b}=\\
\{x:\exists m_x\spa|\ls au,x\rs|^2\leq m_x\ls au,u\rs, \forall u\in\dom{a}\}\oplus
\{y:\exists m_y\spa|\ls bv,y\rs|^2\leq m_y\ls bv,v\rs, \forall v\in\dom{b}\}.
\end{gathered}
\end{equation}
Putting $u=0$ and respectively $v=0$ in \eqref{eq:1}, we see that $$\fdom{a\oplus b}\subseteq \dom{(a_\kr\oplus b_\kr)\sq}.$$ To show $$\fdom{a\oplus b}\supseteq \dom{(a_\kr\oplus b_\kr)\sq},$$ we let $m_{x,y}=2\max(m_x,m_y)$, and use \eqref{eq:1}, \eqref{eq:a} and the convexity of the function $\alpha\mapsto\alpha^2$ on $\mathbb{R}_+$. We have seen consequently that $\fdom{a\oplus b}=\dom{(a_\kr\oplus b_\kr)\sq}$. So the \nam{Krein-von Neumann} extension of $a\oplus b$ exists, and we know that $\fdom{a\oplus b}=\dom{(a\oplus b)_\kr\sq}$. 

To see that $(a\oplus b)_\kr=a_\kr\oplus b_\kr$,
we have to check that $$\dom{(a\oplus b)_\kr\sq}=\dom{(a_\kr\oplus b_\kr)\sq}$$ and furthermore that $$\nm(a_\kr\oplus b_\kr)\sq z\nm^2=\nm (a\oplus b)_\kr\sq z\nm^2$$ holds for all $z\in\dom{(a\oplus b)_\kr\sq}$.

The equality of the domains follows from the above argument.

Now, we prove the required identity. Let $x\oplus y\in\dom{(a\oplus b)_\kr\sq}$.
\begin{equation}\label{eq:2}
\nm(a_\kr\oplus b_\kr)\sq (x\oplus y)\nm^2=\nm (a_\kr\sq \oplus b_\kr\sq)(x\oplus y)\nm^2=\nm a_\kr\sq x\oplus b_\kr\sq y\nm^2=\nm a_\kr\sq x\nm^2+\nm b_\kr\sq y\nm^2
\end{equation}
Now we calculate $\nm(a\oplus b)_\kr\sq(x\oplus y)\nm^2$.
The inequality
\begin{equation}\label{eq:3}\nm(a\oplus b)_\kr\sq(x\oplus y)\nm^2\leq\nm a_\kr\sq x\nm^2+\nm b_\kr\sq y\nm^2
\end{equation} follows immediately from the minimality of the \nam{Krein-von Neumann} extension and the fact that $a_\kr\oplus b_\kr$ is  a positive, selfadjoint extension of $a\oplus b$.

To see the reverse inequality, we consider the following. We can assume that $\|a_\kr\sq x\|^2+\|b_\kr\sq\|^2>0$,therefore we let $$t=\frac{\nm a_\kr\sq x\nm^2}{\nm a_\kr\sq x\nm^2+\nm b_\kr\sq y\nm^2},\quad\mbox{thus}\quad 1-t=\frac{\nm b_\kr\sq y\nm^2}{\nm a_\kr\sq x\nm^2+\nm b_\kr\sq y\nm^2}.$$
\begin{equation*}
\begin{gathered}
\psup_{\begin{ix}\scriptstyle u\in \dom{a},v\in\dom{b}\\\scriptstyle\ls a_\kr u,u\rs+\ls b_\kr v,v\rs\leq 1\end{ix}}\negs\negs\negs{|\ls (a_\kr\sq\oplus b_\kr\sq) (u\oplus v), (a_\kr\sq\oplus b_\kr\sq)(x\oplus y)\rs|^2}\geq\negs\negs\negs\negs
\psup_{\begin{ix}\scriptstyle u\in \dom{a},v\in\dom{b}\\\scriptstyle\ls a_\kr u,u\rs\leq t,\ls b_\kr v,v\rs\leq 1-t\end{ix}}\negs\negs\negs\negs{|\ls a_\kr\sq u, a_\kr\sq x\rs+\ls b_\kr\sq v,b_\kr\sq y\rs|^2}
\end{gathered}
\end{equation*}
Now multiplying $u$ and $v$ by a suitable $\alpha_u,\alpha_v\in\cc$ of absolute value $1$, we continue: 
\begin{equation}
\label{eq:4}
\begin{gathered}
\psup_{\begin{ix}\scriptstyle u\in \dom{a},v\in\dom{b}\\\scriptstyle\ls a_\kr u,u\rs\leq t,\ls b_\kr v,v\rs\leq 1-t\end{ix}}\negs\negs\negs\negs{|\ls a_\kr\sq u, a_\kr\sq x\rs+\ls b_\kr\sq v,b_\kr\sq y\rs|^2}=\negs\negs\negs\negs
\psup_{\begin{ix}\scriptstyle u\in \dom{a},v\in\dom{b}\\\scriptstyle\ls a_\kr u,u\rs\leq t,\ls b_\kr v,v\rs\leq 1-t\end{ix}}\negs\negs\negs\negs{(|\ls a_\kr\sq u, a_\kr\sq x\rs|+|\ls b_\kr\sq v,b_\kr\sq y\rs|)^2}=\\
(\negs\psup_{\begin{ix}\scriptstyle u\in \dom{a},\\\scriptstyle\ls a_\kr u,u\rs\leq t\end{ix}}\negs\negs{|\ls a_\kr\sq u, a_\kr\sq x\rs|}+\negs\negs\psup_{\begin{ix}\scriptstyle v\in\dom{b},\\\scriptstyle\ls b_\kr v,v\rs\leq 1-t\end{ix}}\negs\negs{|\ls b_\kr\sq v, b_\kr\sq y\rs|})^2=\\
t\nm a_\kr\sq x\nm^2+2\sqrt{t(1-t)}\nm a_\kr\sq x\nm\nm b_\kr\sq y\nm+(1-t)\nm b_\kr\sq y\nm^2=
\nm a_\kr\sq x\nm^2+\nm b_\kr\sq y\nm^2
\end{gathered}
\end{equation}
We have used that 
$$\psup_{\begin{ix}\scriptstyle u\in \dom{a},\\\scriptstyle\ls a_\kr u,u\rs\leq t\end{ix}}\negs\negs{|\ls a_\kr\sq u, a_\kr\sq x\rs|^2}=\negs\negs
\psup_{\begin{ix}\scriptstyle u\in \dom{a},\\\scriptstyle\ls a_\kr u,u\rs\leq t\end{ix}}\negs\negs{|\ls a_\kr u, x\rs|^2}=t\nm a_\kr\sq x\nm^2=
\negs\negs\psup_{\begin{ix}\scriptstyle u\in \dom{a_\kr\sq},\\\scriptstyle\ls a_\kr\sq u,a_\kr\sq u\rs\leq t\end{ix}}\negs\negs{|\ls a_\kr\sq u, a_\kr\sq x\rs|^2},$$
and the same for $b_\kr$. Putting together \eqref{eq:2}, \eqref{eq:3} and \eqref{eq:4} we obtain:
$$\nm(a_\kr\oplus b_\kr)\sq (x\oplus y)\nm^2=\nm(a\oplus b)_\kr\sq(x\oplus y)\nm^2$$
completing the proof.
\end{proof}
\end{lem}

\begin{thm}
\label{thm:2}
Let $a$ and $b$ be positive, symmetric operators such that $\fdom{a}\cap\fdom{b}$ is dense in $\hi$, and let $J$ be as in \eqref{eq:0}, then the form sum of $a_\kr$ and $b_\kr$ is $J\aadj J\adj$, i.e. $$a_\kr\fsum b_\kr=J\aadj J\adj.$$ 
\begin{proof}Again we prove that $\dom{(a_\kr\fsum b_\kr)\sq}=\dom (J\aadj J\adj)\sq$, and $ (a_\kr\fsum b_\kr)\sq x=(J\aadj J\adj)\sq x$ for each $x\in\dom{(a_\kr\fsum b_\kr)}$. 

We know that $\dom{(a_\kr\fsum b_\kr)\sq}=\dom{a_\kr\sq}\cap\dom{b_\kr\sq}$, and
$\dom{(J\aadj J\adj)\sq}=\dom J\adj=\fdom{a}\cap\fdom{b}$, as we have seen in the argument following \eqref{eq:0}. Moreover $\dom{a_\kr\sq}=\fdom{a}$ and $\dom{b_\kr\sq}=\fdom{b}$, which implies the desired equality of the domains. 

Using Lemma \ref{thm:1}, we have that
\begin{equation*}
\begin{gathered}
\nm (J\aadj J\adj)\sq x\nm^2=\la J\adj x,J\adj x\ra=\negs\negs\negs
\psup_{\begin{ix}\scriptstyle u\in \dom{a},v\in\dom{b}\\\scriptstyle\la au\oplus bv,au\oplus bv\ra\leq 1\end{ix}}\negs\negs\negs{|\la au \oplus bv,J\adj x\ra|^2}=\negs\negs\negs
\psup_{\begin{ix}\scriptstyle u\in \dom{a},v\in\dom{b}\\\scriptstyle\ls au,u\rs+\ls bv,v\rs\leq 1\end{ix}}\negs\negs\negs{|\ls au+bv,x\rs|^2}=\\
\psup_{\begin{ix}\scriptstyle u\oplus v\in \dom{a}\oplus\dom{b}\\\scriptstyle\ls (a\oplus b)(u\oplus v),u\oplus v\rs\leq 1\end{ix}}\negs\negs\negs{|\ls (a\oplus b)(u\oplus v),x\oplus x\rs|^2}=\\
\nm (a\oplus b)_\kr\sq(x\oplus x)\nm^2=\nm (a_\kr\sq\oplus b_\kr\sq)(x\oplus x)\nm^2=\nm a_\kr\sq x\nm^2+\nm b_\kr\sq x\nm^2.
\end{gathered}
\end{equation*}
Therefore $$\nm (J\aadj J\adj)\sq x\nm^2=\nm a_\kr\sq x\nm^2+\nm b_\kr\sq x\nm^2,$$ which is, by definition, equal to $\nm(a_\kr\fsum b_\kr)\sq x\nm^2$. The theorem is proved.
\end{proof}
\end{thm}

The following theorem is an immediate consequence of Theorem \ref{thm:2}, because for any positive, selfadjoint operator $a$, the \nam{Krein-von Neumann} extension $a_\kr$ and $a$ coincide.
\begin{thm} 
\label{thm:3}
If $a$ and $b$ are positive, selfadjoint operators with $\dom{a\sq}\cap\dom{b\sq}$ dense in $\hi$, then the corresponding operator $J\aadj J\adj$ is just the form sum of $a$ and $b$.
\end{thm}
The previous theorem shows that the following notation is consistent with the notation for the form sum extension. From now on we will use $a\fsum b$ for the above constructed operator $J\aadj J\adj$, even if $a,b$ are positive, symmetric operators. We reformulate Theorem \ref{thm:2} as follows.
\begin{thm}
\label{thm:4}
If $a$ and $b$ are positive, symmetric operators  with $\fdom{a}\cap\fdom{b}$ dense in $\hi$, then $a\fsum b=a_\kr\fsum b_\kr$.
\end{thm}

\begin{rem} Considering the extensions of direct sum of operators, an analogous statement can be proved for the \nam{Friedrichs} extension, as for the \nam{Krein-von Neumann} extension  in Lemma \ref{thm:1}. Namely, if $a,b$ are densely defined, positive, symmetric operators, then 
$$a_\fr\oplus b_\fr=(a\oplus b)_\fr.$$
For the proof we only have to check the equality of the domains of the square root operators.
\begin{equation*}
\begin{split}
\dom(a\oplus b)_\fr\sq=&
\{x\oplus y\in\hi\oplus\hi:\exists x_n\oplus y_n\in\dom{a\oplus b}, x_n\oplus y_n\rightarrow x\oplus y,\\ &\quad\ls (a\oplus b)(x_n\oplus y_n-x_m\oplus y_m),x_n\oplus y_n-x_m\oplus y_m\rs\rightarrow 0\}=\\ &
\{x\oplus y\in\hi\oplus\hi:\exists x_n\in\dom{a}, y_n\in\dom{b}, x_n\rightarrow x,y_n\rightarrow y,\\ &\quad \ls a(x_n-x_m),x_n-x_m\rs+\ls b(y_n-y_m),y_n-y_m\rs\rightarrow 0\}=\\ &
\{x\in\hi:\exists x_n\in\dom{a}, x_n\rightarrow x,\ls a(x_n-x_m),x_n-x_m\rs\rightarrow 0\}\oplus \\ &\oplus\{y\in\hi:\exists y_n\in\dom{b}, y_n\rightarrow y,\ls b(y_n-y_m),y_n-y_m\rs\rightarrow 0\}=\dom{a_\fr\sq}\oplus\dom{b_\fr\sq}
\end{split}
\end{equation*}
\end{rem}

\section{Commutation properties}
In this section we observe that our method constructing the form sum of positive, symmetric operators can preserve some kind of commutation with bounded operators. The ideas used in this section are essentially taken from \cite{sebstoch}, where the commutation property is proved for the \nam{Krein-von Neumann} extension. The situation is as follows: given $E,F\in\bo(\hi)$ and two positive, symmetric operators $a$ and $b$, with $\fdom{a}$ and $\fdom{b}$ dense in $\hi$, such that both $E$ and $F$ leave $\dom{a}$ and $\dom{b}$ invariant. Suppose furthermore that the following equations hold for all $x\in\dom{a}$ and $y\in\dom{b}$:
$$E\adj ax=aF x,\quad F\adj ax=aEx,\quad E\adj by=bF y,\quad F\adj by=bEy.$$
Now, we define $\hat{E}$ and $\hat{F}$ on $\hi_a\oplus\hi_b$ as follows.  
$$\dom{\hat{E}}=\ran{a}\oplus\ran{b},\quad \hat{E}(ax\oplus by)=aEx\oplus bEy,$$and
$$\dom{\hat{F}}=\ran{a}\oplus\ran{b},\quad \hat{F}(ax\oplus by)=aFx\oplus bFy.$$ It is obvious that $\hat{E}$ and $\hat{F}$ leave $\ran{a}\oplus\ran{b}$ invariant. The following lemma shows that both $\hat{E}$ and $\hat{F}$ are well-defined and continuous on a dense subspace of $\hi_a\oplus \hi_b$.
\begin{lem} 
\label{thm:8}
With the notations above, $\hat{E}$ and $\hat{F}$ are well defined, and $\hat{E},\hat{F}\in\bo(\hi_a\oplus\hi_b)$.
\begin{proof} The proof of this lemma could be considerably shortened by referring to the result [Theorem 2 in \cite{sebstoch}]. However, for the sake of completeness we include the detailed proof.
\begin{equation}
\label{eq:5}
\begin{gathered}
\la \hat{F} (ax\oplus by),\hat{F} (ax\oplus by)\ra=\la aFx\oplus bFy,aFx\oplus bFy\ra=\la aFx,aFx\ra+\la bFy,bFy\ra=\\
\ls aFx,Fx\rs+\ls bFy,Fy\rs=\ls E\adj ax,Fx\rs+\ls E\adj by,Fy\rs=\ls ax,EFx\rs+\ls by,EFy\rs=\\
\la ax,aEFx\ra+\la by,bEFy\ra=\la ax\oplus by,aEFx\oplus bEFy\ra\leq\\
\la ax\oplus by,ax\oplus by\ra\sq\la aEFx\oplus bEFy,aEFx\oplus bEFy\ra\sq=\\
\la ax\oplus by,ax\oplus by\ra\sq\la\hat{E}\hat{F}(ax\oplus by),\hat{E}\hat{F}(ax\oplus by)\ra\sq
\end{gathered}
\end{equation}
Substituting $\hat{E}\hat{F}$ for $\hat{F}$, and repeating the argument in \eqref{eq:5}, we obtain 
\begin{equation*}
\la \hat{E}\hat{F} (ax\oplus by),\hat{E}\hat{F} (ax\oplus by)\ra\leq\la ax\oplus by,ax\oplus by\ra\sq\la(\hat{E}\hat{F})^2(ax\oplus by),(\hat{E}\hat{F})^2(ax\oplus by)\ra\sq
\end{equation*}
From this, by induction:
\begin{equation*}
\label{eq:6}
\begin{gathered}
\la \hat{F} (ax\oplus by),\hat{F} (ax\oplus by)\ra\leq
\la ax\oplus by,ax\oplus by\ra^{\frac{1}{2}+\cdots\frac{1}{2^n}}\la(\hat{E}\hat{F})^{\frac{2^n}{2}}(ax\oplus by),(\hat{E}\hat{F})^{\frac{2^n}{2}}(ax\oplus by)\ra^{\frac{1}{2^n}}=\\
\la ax\oplus by,ax\oplus by\ra^{1-\frac{1}{2^n}}\la a(EF)^{\frac{2^n}{2}}x\oplus b(EF)^{\frac{2^n}{2}}y, a(EF)^{\frac{2^n}{2}}x\oplus b(EF)^{\frac{2^n}{2}}y\ra^{\frac{1}{2^n}}=\\
\la ax\oplus by,ax\oplus by\ra^{1-\frac{1}{2^n}}\la (F\adj E\adj)^{\frac{2^n}{2}}ax\oplus (F\adj E\adj)^{\frac{2^n}{2}}by,a(EF)^{\frac{2^n}{2}}x\oplus b(EF)^{\frac{2^n}{2}}y\ra^{\frac{1}{2^n}}=\\
\la ax\oplus by,ax\oplus by\ra^{1-\frac{1}{2^n}}\ls ax\oplus by, (EF)^{2^n}x\oplus (EF)^{2^n}y\rs^{\frac{1}{2^n}}\leq\\
\la ax\oplus by,ax\oplus by\ra^{1-\frac{1}{2^n}}\nm ax\oplus by\nm^{\frac{1}{2^n}}\nm(EF)^{2^n}x\oplus (EF)^{2^n}y\nm^{\frac{1}{2^n}}=\\
\la ax\oplus by,ax\oplus by\ra^{1-\frac{1}{2^n}}\nm ax\oplus by\nm^{\frac{1}{2^n}}\nm((EF)^{2^n}\oplus (EF)^{2^n})(x\oplus y)\nm^{\frac{1}{2^n}}\leq\\
\la ax\oplus by,ax\oplus by\ra^{1-\frac{1}{2^n}}\nm ax\oplus by\nm^{\frac{1}{2^n}}\nm(EF)^{2^n}\oplus (EF)^{2^n}\nm^{\frac{1}{2^n}} \nm x\oplus y\nm^{\frac{1}{2^n}}=\\
\la ax\oplus by,ax\oplus by\ra^{1-\frac{1}{2^n}}\nm ax\oplus by\nm^{\frac{1}{2^n}}\nm(EF\oplus EF)^{2^n}\nm^{\frac{1}{2^n}} \nm x\oplus y\nm^{\frac{1}{2^n}}
\end{gathered}
\end{equation*}
If we take the limit $n\rightarrow\infty$, we obtain:
$$\la \hat{F} (ax\oplus by),\hat{F} (ax\oplus by)\ra\leq r(EF\oplus EF)\la ax\oplus by,ax\oplus by\ra,$$
where $r(EF\oplus EF)$ stands for the spectral radius of $EF\oplus EF$. And this is enough to prove both statements for $\hat{F}$. The proposition for $\hat{E}$ can be proved analogously. (To be very precise, we have shown that $\hat{E}$ and $\hat{F}$ are continuously defined on a dense subspace of $\hi_a\oplus\hi_b$, but they are automatically extended to the whole space.)
\end{proof}
\end{lem}
Now, we compute the adjoints of $\hat{E}$ and $\hat{F}$ in $\bo(\hi_a\oplus\hi_b)$:
\begin{lem} \label{thm:9}
$\hat{E}\adj=\hat{F}$ and $\hat{F}\adj=\hat{E}$.
\begin{proof}
It is enough to prove  $\hat{F}\adj=\hat{E}$, as  $\hat{E},\hat{F}\in\bo(\hi_a\oplus\hi_b)$. We check that $\hat{F}\adj x=\hat{E}x$ on the dense subspace $\ran{a}\oplus\ran{b}$. Let $ax\oplus by\in\ran{a}\oplus\ran{b}$, then for all  $au\oplus bv\in\ran{a}\oplus\ran{b}$
\begin{equation*}
\begin{split}
\la au\oplus bv,\hat{F}\adj(ax\oplus by)\ra=\la \hat{F}(au\oplus bv),ax\oplus by\ra=\la aFu\oplus bFv,ax\oplus by\ra=\\
\la aFu,ax\ra+\la bFv,by\ra=\ls aFu,x\rs+\ls bFv,y\rs=\ls E\adj au,x\rs+\ls E\adj bv,y\rs=\ls au,Ex\rs+\ls bv,Ey\rs=\\
\la au,aEx\ra+\la bv,bEy\ra=\la au\oplus bv,aEx\oplus bEy\ra=\la au\oplus bv,\hat{E}(ax\oplus by)\ra, 
\end{split}
\end{equation*}
and that was to be proved.
\end{proof}
\end{lem}

\begin{thm}
\label{thm:5}
 Let $a,b$ be positive, symmetric operators with $\fdom{a}\cap\fdom{b}$ dense in $\hi$, and suppose  that $E,F\in\bo(\hi)$, such that both $E$ and $F$ leave $\dom{a}$ and $\dom{b}$ invariant, and for all $x\in\dom{a}$ and $y\in\dom{b}$
$$E\adj ax=aF x,\quad F\adj ax=aEx,\quad E\adj by=bF y,\quad F\adj by=bEy.$$
Then $$E\adj(a\fsum b)\subseteq(a\fsum b)F\quad\mbox{and}\quad F\adj(a\fsum b)\subseteq (a\fsum b)E.$$
\begin{proof}
First we show the following:
$$E\adj J\subseteq J\hat{F},\quad F\adj J\subseteq J\hat{E},\quad \hat{E} J\adj\subseteq J\adj E,\quad\hat{F} J\adj\subseteq J\adj F.$$
Indeed, let $ax\oplus by\in\ran{a}\oplus\ran{b}$, then
$$J\hat{F}(ax\oplus by)=J(aFx\oplus bFy)=aFx+bFy=E\adj ax+E\adj by=E\adj(ax+by)=E\adj J(ax\oplus by).$$
Observing the domains, we have consequently $E\adj J\subseteq J\hat{F}$. An analogous proof can be given for $F\adj J\subseteq J\hat{E}$. For the remaining inclusions, we write:
$$\hat{E}J\adj=\hat{F}\adj J\adj \subseteq (J\hat{F})\adj\subseteq (E\adj J)\adj=J\adj E,$$
as $E$ is bounded, hence $\hat{E}J\adj\subseteq J\adj E$, and with the same reasoning $\hat{F}J\adj\subseteq J\adj F$.

Finally we turn to the proof of the theorem. Using the previously proved statement, we have
$$E\adj J\aadj\subseteq (J\adj E)\adj\subseteq(\hat{E}J\adj)\adj=J\aadj \hat{E}\adj=J\aadj \hat{F}.$$
Note that we have used that $\hat{E}$ is continuous according to Lemma \ref{thm:8}.
We complete the proof by writing
$$E\adj(a\fsum b)=E\adj J\aadj J\adj\subseteq J\aadj\hat{F} J\adj\subseteq J\aadj J\adj F=(a\fsum b)F,$$ that is $E\adj(a\fsum b)\subseteq(a\fsum b)F$, and with the same argument $F\adj(a\fsum b)\subseteq(a\fsum b)E$.
\end{proof}
\end{thm}

The following result, which is just a special case of Theorem \ref{thm:5} with $E=F=S=S\adj$, shows the reason why we talk about ``commutation properties'' above.
\begin{thm}Let $S$ be a bounded, selfadjoint operator over the Hilbert space $\hi$, such that $S$ leaves both $\dom{a}$ and $\dom{b}$ invariant, and furthermore
$$Sax=aS x,\quad S by=bSy$$
hold for all $x\in\dom{a}$ and $y\in\dom{b}$. Also, assume that $\fdom{a}\cap \fdom{b}$ is dense in $\hi$. Then
$$S(a\fsum b)\subseteq(a\fsum b)S.$$
\end{thm}

In Theorem \ref{thm:5}, we required that the bounded operators $E,F$ leave some subspaces invariant. In some cases, we might not know that such ``big'' subspaces are invariant, perhaps because they are not invariant at all, but we may find smaller subspaces whose invariance can be checked. We try to handle this problem, with the following theorem. 
\begin{thm} 
\label{thm:6}
Let $a$ and $b$ be positive, symmetric operators with $\fdom{a}\cap\fdom{b}$ dense in $\hi$, and suppose that $D\subseteq\dom{a}\cap\dom{b}$ is a linear manifold. Then $a\rest_D\fsum b\rest_D=a\fsum b$ if and only if for all $x\in\hi$
\begin{equation}\label{eq:7}
\psup_{\begin{ix}\scriptstyle u\in \dom{a},\\\scriptstyle\ls au,u\rs\leq 1\end{ix}}\negs{|\ls au,x\rs|^2}+
\psup_{\begin{ix}\scriptstyle v\in \dom{b},\\\scriptstyle\ls bv,v\rs\leq 1\end{ix}}\negs{|\ls bv,x\rs|^2}=
\negs\psup_{\begin{ix}\scriptstyle u\in {D},\\\scriptstyle\ls au,u\rs\leq 1\end{ix}}\negs{|\ls au,x\rs|^2}+
\negs\psup_{\begin{ix}\scriptstyle v\in {D},\\\scriptstyle\ls bv,v\rs\leq 1\end{ix}}\negs{|\ls bv,x\rs|^2}
\end{equation}
\begin{proof} Before all, observe that $\fdom{a}\subseteq\fdom{a\rest_D}, \fdom{b}\subseteq\fdom{b\rest_D}$, indeed:
\begin{equation}
\label{eq:8}
\begin{gathered}
\fdom{a}=\{y\in\hi:\exists m_y |\ls ax,y\rs|^2\leq m_y\ls ax,x\rs,\forall x\in \dom{a}\}\subseteq\\
\{y\in\hi:\exists m_y |\ls ax,y\rs|^2\leq m_y\ls ax,x\rs,\forall x\in D\}=\fdom{a\rest_D},
\end{gathered}
\end{equation}
and the same for $\fdom{b}$ and $\fdom{b\rest_D}$.

Suppose now that condition \eqref{eq:7} is satisfied. Then for the reverse inclusion  $\fdom{a}\cap\fdom{b}\supseteq\fdom{a\rest_D}\cap\fdom{b\rest_D}$ we let $x\in \fdom{a\rest_D}\cap\fdom{b\rest_D}$, which is the same as saying that the right hand side of \eqref{eq:7} is finite for this $x$. But then, from assumption \eqref{eq:7} it follows that the left hand side of \eqref{eq:7} is also finite, implying $x\in \fdom{a}\cap\fdom{b}$.
 By our construction for the form sum $$\dom((a\fsum b)\sq)=\fdom{a}\cap\fdom{b},\quad\mbox{and}\quad \dom(a\rest_D\fsum b\rest_D)\sq=\fdom{a\rest_D}\cap\fdom{b\rest_D},$$
hence $\dom(a\fsum b)\sq=\dom(a\rest_D\fsum b\rest_D)\sq$. Let $x\in \dom(a\fsum b)\sq$, then by the proof of Theorem  \ref{thm:2} and \eqref{eq:-0}
\begin{equation}
\begin{gathered}
\label{eq:9}
\nm (a\fsum b)\sq x\nm^2=\nm a_\kr\sq x\nm^2+\nm b_\kr\sq x\nm^2=\negs
\psup_{\begin{ix}\scriptstyle u\in \dom{a},\\\scriptstyle\ls au,u\rs\leq 1\end{ix}}\negs{|\ls au,x\rs|^2}+\negs\psup_{\begin{ix}\scriptstyle v\in \dom{b},\\\scriptstyle\ls bv,v\rs\leq 1\end{ix}}\negs{|\ls bv,x\rs|^2}=\\
\psup_{\begin{ix}\scriptstyle u\in D,\\\scriptstyle\ls au,u\rs\leq 1\end{ix}}\negs{|\ls au,x\rs|^2}+\negs\psup_{\begin{ix}\scriptstyle v\in D,\\\scriptstyle\ls bv,v\rs\leq 1\end{ix}}\negs{|\ls bv,x\rs|^2}=\nm (a\rest_D)_\kr\sq x\nm^2+\nm (b\rest_D)_\kr\sq x\nm^2=\nm (a\rest_D\fsum b\rest_D)\sq x\nm^2.
\end{gathered}
\end{equation}
Consequently we have $a\rest_D\fsum b\rest_D=a\fsum b$.

For the reverse direction, we suppose that $a\rest_D\fsum b\rest_D=a\fsum b$. Then for all $x\in\dom(a\fsum b)\sq$ $\nm (a\fsum b)\sq x\nm^2=\nm (a\rest_D\fsum b\rest_D)\sq x\nm^2$, and the same argument as in \eqref{eq:9} shows that \eqref{eq:7} is satisfied.
\end{proof}
\end{thm} 
\section{Further results and remarks}
Our construction for the form sum is based on the idea used when constructing the \nam{Krein-von Neumann} extension $a_\kr$ of a positive, symmetric operator $a$. Analogously we consider the following situation. We suppose that $\dom{a}$ and $\dom{b}$ are dense. Again we have the Hilbert space $\hi_a\oplus\hi_b$, and we define analogously as in  \cite{sebprok2}, \cite{sebprok1} $$Q:\hi\rightarrow\hi_a\oplus\hi_b,\mbox{ with } \dom{Q}=\dom{a}\cap\dom{b},\enskip Qx=ax\oplus bx.$$ Obviously $Q$ is a restriction of $J\adj$. The question is, what can be said about $Q\adj Q\aadj$.
\begin{thm} 
\label{thm:7}
Suppose that $a$ and $b$ are positive, symmetric operators, and $\dom{a}\cap\dom{b}$ is dense in $\hi$. Then $Q\adj Q\aadj=(a+b)_\fr$.
\begin{proof}
First we show that under these circumstances $Q\adj Q\aadj$ exists and is a positive, selfadjoint operator. From the \nam{von Neumann} theorem, it is clear that if $Q\adj Q\aadj$ exists then it is selfadjoint, and obviously positive. $Q\adj$ exists, since $\dom{Q}$ is dense. We compute $\dom{Q\adj}$, and as it will be dense, we conclude that $Q\aadj$ exists. 
First we compute $Q\adj$ on $\ran{a}\oplus\ran{b}$.
Let $ax\oplus by\in\ran{a}\oplus\ran{b}$ and $z\in\dom{a}\cap\dom{b}$
\begin{equation*}
\begin{gathered}
\la Qz,ax\oplus by\ra=\la az\oplus bz,ax\oplus by\ra=\la az,ax\ra+\la bz,by\ra=\ls az,x\rs+\ls bz,y\rs=\\
\ls z,ax\rs+\ls z,by\rs=\ls z,ax+by\rs,
\end{gathered}
\end{equation*}
which shows that $\ran{a}\oplus\ran{b}\subseteq \dom{Q\adj}$ and $Q\adj(ax\oplus by)=ax+by$. Therefore $Q\adj$ is densely defined. We see that $Q\adj Q\aadj$ is an extension of $a+b$:  
$$Q\adj Q\aadj z=Q\adj Q z=Q\adj(az\oplus bz)=az+bz.$$

Because of the extremality of the \nam{Friedrichs} extension, we only have to prove that $$\dom(a+b)_\fr\sq=\dom (Q\adj Q\aadj)\sq.$$
We can write
\begin{equation*}
\begin{gathered}
\dom{(Q\adj Q\aadj)}\sq=\dom{Q\aadj}=\dom{\bar{Q}}=\{y\in\hi:\exists y_n\in\dom{Q},y_n\rightarrow y,Qy_n\mbox{ convergent}\}=\\
\{y\in\hi:\exists y_n\in\dom{Q},y_n\rightarrow y, \la ay_n\oplus by_n-ay_m\oplus by_m,ay_n\oplus by_n-ay_m\oplus by_m\ra\rightarrow 0\}=\\
\{y\in\hi:\exists y_n\in\dom{a}\cap\dom{b},y_n\rightarrow y, \ls a(y_n-y_m),y_n-y_m\rs +\ls b(y_n-y_m),y_n-y_m\rs\rightarrow 0\}=\\
\{y\in\hi:\exists y_n\in\dom(a+b),y_n\rightarrow y, \ls (a+b)(y_n-y_m),y_n-y_m\rs\rightarrow 0\}=\dom(a+b)_\fr\sq,
\end{gathered}
\end{equation*}
which remained to complete the proof.
\end{proof}
\end{thm}

Finally, we examine the connection between different extensions of the operator sum. Supose that $A$ and $B$ are positive, selfadjoint operators, and let $A+B$ denote the operator sum on $D=\dom{A}\cap\dom{B}$. Suppose that $D$ is dense in $\hi$, so that the \nam{Friedrichs} extension $(A+B)_\fr$ of $A+B$ exists. \nam{Kato} \cite{kato} shows an example when $A\fsum B\neq (A+B)_\fr$. Analogously, one can examine the connection between $A\fsum B$ and $(A+B)_\kr$. We will prove that in general $A\fsum B\neq (A+B)_\kr$. Note if we assume only that $\dom{A\sq}\cap\dom{B\sq}$ is dense in $\hi$ -- assuring the existence of $A\fsum B$ -- the \nam{Krein-von Neumann} extension will still exist. Indeed, it is easy to see that 
\begin{gather*}
\fdom{A+B}=\{y\in\hi:\exists m_y\: |\ls (a+b)x,y\rs|^2\leq m_y\ls(a+b)x,x\rs,\ \forall x\in D\}\supseteq\\ \fdom{A}\cap\fdom{B}=\dom{A\sq}\cap\dom{B\sq},
\end{gather*}
so $\fdom{A+B}$ is dense in $\hi$. However, it may happen that $\dom{A\sq}\cap\dom{B\sq}$ is dense in $\hi$ while $\dom{A}\cap\dom{B}=\{\nul\}$. In this case $A\fsum B\neq (A+B)_\kr=\nul$, providing a trivial counter-example. For this reason, in the sequel we keep the assumption that $D$ is dense in $\hi$.
\begin{exa}
Consider the following example. Let $a$ be a densely defined, closed, symmetric operator with positive lower bound. Suppose moreover that $a$ is not selfadjoint. Then the deficiency index $\dim(\ker a\adj)$ of $a$ is greater than zero. Also, there are infinitely many selfadjoint extensions of $a$, which are restrictions of $a\adj$. Among these the \nam{Friedrichs} extension $a_\fr$ is the largest, and the \nam{Krein-von Neumann} $a_\kr$  is the smallest one with respect to the usual ordering of positive, selfadjoint operators. Consider $a_\kr$ and $a_\fr$, both are positive and selfadjoint, and $D=\dom{a_\kr}\cap\dom{a_\fr}\supseteq \dom{a}$, so $D$ is dense in $\hi$. Furthermore, we have that $a_\kr\fsum a_\fr=2a_\fr$,
because $$\dom{a_\kr\sq}\cap\dom{a_\fr\sq}=\dom{a_\fr\sq},\mbox{ and }\nm a_\kr\sq x\nm^2=\nm a_\fr\sq x\nm^2$$ for all $x\in\dom{a_\fr\sq}$. On the other hand, $(a_\kr+a_\fr)_\kr=2a_\kr$, because $a_\kr+a_\fr$ is a symmetric extension of $2a$, hence $2a_\kr=(2a)_\kr\leq(a_\kr+a_\fr)_\kr$. Conversely, $(a_\kr+a_\fr)_\kr\leq 2a_\kr$, because $a_\kr+a_\fr$ is a restriction of $2a_\kr$. Thus we have that $a_\kr\fsum a_\fr\neq(a_\kr+a_\fr)_\kr$, as desired.
\end{exa}
\begin{exa}
A similar approach can provide an example when $A\fsum B\neq (A+B)_\fr$. The example above fails as $a_\kr\fsum a_\fr=2 a_\fr$ and  $(a_\kr+a_\fr)_\fr=2 a_\fr$ as well. However, take any intermediate extension $a_\im$ of $a$ instead of $a_\fr$. Then we have $a_\kr\fsum a_\im\leq 2a_\im$ because $$\dom{(a_\kr\fsum a_\im)\sq}=\dom{a_\kr\sq}\cap\dom{a_\im\sq}=\dom{a_\im\sq}$$
and 
$$\nm(a_\kr\fsum a_\im)\sq x\nm^2=\nm a_\kr\sq x\nm^2+\nm a_\im\sq x\nm^2\leq\nm(2a_\im)\sq x\nm^2$$
for all $x\in\dom{a_\im\sq}$. Furthermore, $(a_\kr+a_\im)_\fr\geq 2a_\im$ because both $a_\kr+a_\im$ and $2a_\im$ are extensions of $2a_\im\rest_{\dom{a_\kr}\cap\dom{a_\im}}=a_\kr+a_\im$, here we have used that 
$$a_\kr\rest_{\dom{a_\kr}\cap\dom{a_\im}}=a\adj\rest_{\dom{a_\kr}\cap\dom{a_\im}}=a_\im\rest_{\dom{a_\kr}\cap\dom{a_\im}},$$
so the inequality follows from the extremality of the \nam{Friedrichs} extension. Thus we have
\begin{equation}
\label{eq:11}
a_\kr\fsum a_\im\leq 2a_\im\leq (a_\kr+a_\im)_\fr.
\end{equation}
 How can we assure that equality does not hold at both inequalities in \eqref{eq:11}? It is easy to see from the above that a sufficient condition for $a_\im$ is that the form $q_{a_\im}$ of $a_\im$ is not a restriction of the form $q_{a_\kr}$ of $a_\kr$. When $\dim{(\ker a\adj)}>0$, such an $a_\im$ is always available (see \cite{alsim}). Just take any strictly positive, closed form $q_0$ on $\ker a\adj$ (e.g. the original inner product) and define a new form $q$ on $\ker a\adj+\dom{a_\fr\sq}$ as follows
$$q(x+y)=q_0(x)+\nm a_\fr\sq y\nm^2,\quad x\in\ker a\adj, y\in\dom{a_\fr\sq}.$$
We have used that  $\ker{a\adj}\cap\dom{a_\fr\sq}=\{\nul\}$. Using the representation theorem, we get the required $a_\im$. (Note that $a_\kr$ belongs to the choice $q_0\equiv 0$.) Thus we see that a desired counter-example can be given whenever $\dim\ker a\adj>0$.
\end{exa}
\nocite{glimmjaffe}

\end{document}